<div class="moz-text-flowed" style="font-family: -moz-fixed">
\input amstex
\input epsf
\documentstyle{amsppt}
\magnification=1200
\NoBlackBoxes
\nologo
\tolerance 1000
\hsize 16.4 truecm
\vsize 22.9 truecm

\global\topskip=0pt
\baselineskip 18pt

\define \Ga {\Gamma}
\define \la {\lambda}
\define \pa {\partial}
\define \bR {\Bbb R}
\define \bC {\Bbb C}
\define \bF {\Bbb F}
\define \al {\alpha}

\def \ga {\gamma}
\define \Si {\Sigma}
\define \si {\sigma}
\def \bP {\Bbb RP}
\def \bcP {\Bbb CP}
\define \gkm {G_{k,m}}
\def \gk1 {G_{k,n+1}}

\define \gtm {G_{2,m}}
\define \gtn {G_{2,n}}
\define \gkn {G_{k,n}}
\def \gk1 {G_{k,n+1}}

\vskip 1.5cm

\topmatter

\title
On two conjectures concerning convex curves
\endtitle

\author
       V.~Sedykh  and B.~Shapiro
\endauthor


\affil
       Department of Mathematics, University of Oil and Gas (Gubkin) \\
Moscow 117036, Russia, {\tt sedykh\@mccme.ru}
       Department of Mathematics, University of Stockholm\\
S-10691, Sweden, {\tt shapiro\@math.su.se}
\endaffil


\abstract  In this paper we recall two basic conjectures on the
developables of convex projective curves, prove one of them and
disprove the other in the first nontrivial case of curves in $\bP^3$.
Namely, we show that  i) the tangent developable of any convex curve
in $\bP^3$ has degree $4$ and ii) construct an example of $4$ tangent
lines to a convex curve in $\bP^3$ such that no real line intersects
all four of them.
The question (discussed in \cite {EG1} and \cite {So4})
whether the second conjecture is true in the special case
of  rational normal curves still remains open.

\endabstract
\rightheadtext{Two conjectures }
\leftheadtext{  V.~Sedykh and B.~Shapiro }
\endtopmatter
\document
\heading {\S 1. Introduction and results} \endheading

We start with some  important notions.
\medskip
{\smc Main definition.} A smooth closed curve $\ga: S^1\to \bP^n$ is
called {\it locally convex} if the local multiplicity of intersection
of $\ga$ with any hyperplane $H\subset \bP^n$ at any of the intersection
points does not exceed $n=\dim \bP^n$ and {\it globally convex} or
just {\it convex} if the above condition holds for the global multiplicity,
i.e for the sum of local multiplicities.
\medskip
  Local convexity of $\ga$ is a simple requirement  of
nondegeneracy of the
osculating Frenet $n$-frame of $\ga$, i.e. of the linear independence
of $\ga^\prime(t),...,\ga^{(n)}(t)$ for any $t\in S^1$. Global convexity
is a nontrivial property equivalent to the fact that the
$(n+1)$-tuple of $\ga$'s homogeneous coordinates forms a
Tschebychev system of functions, see e.g. \cite {KS}.
  The simplest examples of convex curves are the rational
normal curve $\rho_{n}: t\mapsto (t, t^2, \ldots, t^n)$ (in some
affine coordinates on $\bP^n$) and the standard trigonometric curve
$\tau_{2k}: t\mapsto (\sin t, \cos t, \sin 2t, \cos 2t,\ldots, \sin
2kt, \cos 2kt)$ (in some affine coordinates on $\bP^{2k}$).
\medskip
{\smc Definition.} The {\it $k$-th  developable} $D_{k}(\ga)$
of a curve
$\ga: S^1\to \bP^n$ is the union of all $k$-dimensional osculating
subspaces to $\ga$. The hypersurface $D(\ga)=D_{n-2}(\ga)$ is called
the {\it developable hypersurface} of $\ga$. ($D(\ga)$ was called the
{\it discriminant} of $\ga$ in \cite {SS}.)
\medskip
Note that  $D_{k}(\ga)$  can be considered as the image of the
natural associated
map $\ga_{k}: S^1\times \bP^k \to \bP^n$.
\medskip
{\smc Definition.}
     Let $M$ be a compact manifold of some dimension $l\le n$ and $\phi:
     M\to \bP^n$ be a smooth map. By the {\it degree} of $\phi(M)$ we 
understand the
     supremum of the number of its intersection points with generic
     $(n-l)$-dimensional subspaces. Recall that an
     $(n-l)$-dimensional subspace $L$ is called {\it generic} w.r.t 
$\phi(M)$ if the
     intersection $L\cap \phi(M)$ consists of a finite number of
     points and at each such point the tangent spaces to $\phi(M)$ and
     $L$ are transversal, i.e.  their sum coincides with the
     whole tangent space to $\bP^n$ at this point. Notice that the
     degree of $\phi(M)$ can be infinite. If there are no
     generic $(n-l)$-dimensional subspaces for $\phi(M)$ (for example,
     if $\dim \phi(M) <\dim (M)$) then we set the degree of $\phi(M)$
     equal to zero. (It is rather obvious that if the Jacobian of
     $\phi$ is nondegenerate at least at one point of $M$ then generic
     $(n-l)$-dimensional subspaces exist.)
     \medskip
     In particular, one can consider the degree of $D_{k}(\ga)$
     which is positive unless $\dim D_{k}(\ga) < k+1.$
\medskip
{\smc Remark.} By definition $\ga$ is convex if and only if
the degree of $\ga=D_{0} (\ga)$ (considered as the image of the map
$\ga: S^1\to \bP^n$) equals $n$. It is well-known that
$\ga$ is convex if and only if the dual curve $\ga^*\in (\bP^n)^{*}$
is convex, see e.g. \cite {A1}. Thus  $\ga$ is convex if
and only if the degree of $D_{n-1}(\ga)$ as the image of associated
map $\ga_{n-1}: S^{1}\times \bP^{n-1} \to \bP^n$ equals $n$. The
following conjecture generalizes the latter fact.

\medskip
{\bf  Conjecture on the $k$-th developable.} A curve $\ga: S^1\to
\bP^n$ is convex
if and only if for some (and then for all) $k=0, \ldots , n-1$ the
degree of its $k$-th developable $D_{k}(\ga)$ equals  $(k+1)(n-k)$.
\medskip
(Note that $(k+1)(n-k)$ is the dimension of
the Grassmannian of projective $k$-dimensional subspaces in $P^n$.)

\medskip

As a special case one gets
\medskip
{\bf   Conjecture on the developable hypersurface.} A curve $\ga: S^1\to \bP^n$
is convex  if and only if the degree of its developable hypersurface
$D(\ga)\subset \bP^n$ equals $2n-2$.

\medskip
{\smc Remark.} The conjecture on the $k$-th developable can be verified  by a
straightforward degree count in the case when $\ga$ is the rational 
normal curve.
Our first result is the following

\medskip
{\bf Theorem A.} A curve $\ga:S^1\to \bP^3$ is convex if and only if
the degree of its developable hypersurface $D(\ga)$ equals $4$.
\medskip
In order to formulate the second conjecture recall that
     for  a  generic $(k+1)(n-k)$-tuple of $k$-dimensional complex subspaces in
     $\bcP^n$ there exist  $\sharp_{k,n}=\frac {1!2!\ldots
(n-k-1)!((k+1)(n-k))!}{(k+1)!(k+2)!\ldots(n)!}$ projective complex
     subspaces of dimension $(n-k-1)$ in $\bcP^n$ intersecting each of
the above $k$-dimensional susbspaces. This is a classical result due to
H.~Schubert, see \cite {Sch}.
(The number $\sharp_{k,n}$ is the degree of the Grassmanian of
projective $k$-dimensional subspaces in $P^n$
considered as a projective variety embedded using
Pl\"ucker coordinates.)

\medskip

{\bf Conjecture on total reality. } For the real rational normal curve
$\rho_{n} : S^1\to
\bP^n$ and any  $(k+1)(n-k)$-tuple of pairwise distinct {\it real} projective
$k$-dimensional  osculating subspaces to $\rho_{n}$ there exist
$\sharp_{k,n}$ {\it real} projective
     subspaces of dimension $(n-k-1)$ in $\bP^n$ intersecting each of the
     above osculating subspaces.
\medskip
The latter conjecture is sometimes referred to as the B.\&M.Shapiro
conjecture and has the
following appealing interpretation in terms of real algebraic geometry.
Consider a generic degree $n$ rational curve: $\mu: \bcP^1\to \bcP^k$.
Such a curve has exactly $(k+1)(n-k)$ inflection points. (An
{\it inflection point} of a space curve $\mu$ is a point where the osculating
Frenet frame $(\mu',\mu'',\ldots,\mu^{(k)})$ is degenerate.)
\medskip
{\bf  Conjecture on total reality restated. }  In the above
notation assume additionally that the inverse images of
all inflection points lie on $\bP^1\subset \bcP^1$. Then $\mu$ is a
{\it real} rational curve up to a projective transformation in
$\bcP^k$, i.e. there exists a projective transformation of $\bcP^k$
making $\mu$ into a curve invariant under the complex conjugations in
the preimage and the image.
\medskip
The conjecture on total reality in the case
$k=1$, i.e. for rational functions in one variable
was  settled in \cite {EG1}. Extensive numerical support for
its validity  for $k > 1$ can be found in  \cite {So4} and its
intriguing generalization to the case of manifolds of (in)complete
flags is contained in \cite {RSSS}.
\medskip
Earlier the second author proposed a {\bf generalized conjecture on total
reality}
saying that the above conjecture should be valid for
all convex curves in $\bP^n$ (and not just for the rational normal
curve), see \cite {So4}. But this generalized
conjecture on total reality  fails as shown by the next result.
\medskip
{\bf Theorem B.} There exists a convex curve $\ga:S^1\to \bP^3$ and
a $4$-tuple of its tangent lines with no real lines in
$\bP^3$ intersecting all four tangents.
\medskip
{\it The most delicate open problem related to the present paper is to
find a natural subset  of the set of convex curves for which the conjecture of
total reality might hold.}

\medskip
{\smc Remarks.} One of the implications in Theorem A was earlier
proven in \cite {SS}. The conjecture on total reality
was discussed in a number of papers, see \cite {EG1}, \cite
{So1}-\cite {So4}. The structure of
the paper is as follows. \S 2-3 contain the proofs of Theorems A and B
resp. \S 4 discusses the above conjectures in terms of Schubert
calculus. In particular, it contains a stronger version of the total
reality conjecture closely related to the transversality conjecture due
to F.~Sottile.
\medskip
{\smc Acknowledgements.} The first author is sincerely grateful to
the Department of Mathematics, University of Stockholm for the
hospitality and to the Department of Mathematics, Royal Institute of
Technology for the financial support of his visit to Sweden in April
2002 when this project was started. The second author wants to
emphasize the importance of his contacts with A.~Gabrielov,
A.~Eremenko, M.~Shapiro and, especially, with F.~Sottile.

\heading {\S 2. Proof of Theorem A}\endheading

    The idea of the proof of Theorem A is as follows. We will show
    that the degree of the developable hypersurface is either
    identically equal to $4$ or strictly greater than $4$ for all curves within
    any connected component of the space of locally convex curves in $\bP^3$.
    (The main technical tool is Proposition 1 below.)
     After that in order to settle Theorem A one just has to calculate the
     degree of the developable hypersurface for
     some suitable representative in each such connected component, see
     Proposition 2.

\medskip

Take a smooth $1$-parameter family
  of smooth closed curves $\ga_t: S^1\to \bP^3,\; t \in
[0,1]$ and assume that:

\medskip
($\star$) for any $t\in (0,1]$ the curve $\ga_t$ is locally convex 
and the degree of its
developable hypersurface $D(\ga_t)$ equals $4$;

($\star\star$) the degree of the developable hypersurface
$D_0=D(\ga_0)$  exceeds $4$.

\medskip
The next technical statement is the key step in the proof of Theorem
A.
\medskip
{\smc Proposition 1.} {\it Under the above assumptions the curve $\ga_0$
is not locally convex, i.e. has inflections.}

\medskip
To prove this result  we assume that  $\ga_0$ is locally
convex and get a contradiction with ($\star$) by finding a line
intersecting $D_{0}$ transversally more than $4$ times.

\medskip
First of all note that any curve $\tilde\ga$ sufficiently close to
$\ga_0$ in $C^1$-topology is locally convex. Hence, every irreducible 
component of
a germ of  the
developable hypersurface of  $\tilde\ga$ at any point   is either
a smooth germ  or is diffeomorphic
to a germ of a half-cubical cuspidal edge with singularities at
$\tilde\ga$. The number of local irreducible components of a germ of 
the developable
hypersurface at any point is finite. Moreover, we can assume that 
that the family $\ga_t$
is generic in the sense that the front  of tangent planes  (i.e. the
developable hypersurface of the dual curve) to the
curve $\ga_t$ has only generic perestroikas at the moment $t=0$.
These perestroikas were listed in \cite{A2}, p.2.

The next lemma collects  restrictions valid for a locally convex 
$\ga_{0}$ satisfying
($\star$) and ($\star\star$) in the increasing order of strength. 
(Recall that our aim is to prove that
such a $\ga_{0}$ cannot exist.)

\medskip
{\smc 2.1. Lemma.} {\it If $\ga_0$ is locally convex and satisfies
($\star$) and ($\star\star$) then

1)  $\ga_{0}$ cannot  intersect any of its osculating planes
transversally;

2) $\ga_0$ has no self-intersections;

3)  $\ga_0$ has exactly one common point
with any of its tangent lines;

4)  $\ga_0$ has exactly one common point with any of its osculating
planes.}

Obviously, $4)\Rightarrow 3) \Rightarrow 2) \Rightarrow 1).$

\demo{Proof} To prove 1) notice that if $\ga_{0}$ transversally 
intersects the plane
osculating it at a point $p$ at some point $q\neq p$,
then for any sufficiently small $t$ the  plane osculating the curve
$\ga_t$ at some point $\tilde p$ close to $p$ transversally intersects
$\ga_t$ at some point $\tilde q$ close to $q$. By continuity the point
$\ga_t(\tilde q)$ does not lie on the tangent line to $\ga_t$ at
$\tilde p$. The line passing through the points
$\ga_t(\tilde p),\ga_t(\tilde q)$ intersects $D(\ga_t)$ with the
multiplicity $>4$ contradicting ($\star$).

To prove 2) notice that the osculating plane to one of the local 
irreducible components
of the curve $\ga_0$ at its self-intersection point intersects another
irreducible component transversally (since our family $\ga_t$ is
chosen generic by assumption). But this is impossible by the above item 1).

To prove 3) assume that the tangent line to $\ga_0$ at a point $p$
intersects $\ga_0$ at some point $q\neq p$. Then
$\ga_0(p)\neq \ga_0(q)$ by item 2)  and
the  plane osculating $\ga_0$ at $p$ is tangent to $\ga_0$ at $q$
by item 1). The  plane osculating $\ga_0$ at $q$ is transversal to
the line passing through the points $\ga_0(p)$ and $\ga_0(q)$ (since
the family $\ga_t$ is generic). Hence for any sufficiently small $t$
there exists a point $\tilde p$ close to $p$ and a line $l$ in $\bP^3$
passing through $\ga_t(\tilde p)$ such that: a) $l$ belongs to the
plane osculating $\ga_t$ at $\tilde p$; b) $l$ is not tangent to 
$\ga_t$ at $\tilde p$;
c) $l$ intersects the irreducible component of the surface $D(\ga_t)$ at the
point $\ga_t(\tilde q)$ at two smooth points. The total multiplicity
of the intersection $l\cap D(\ga_t)$ exceeds $4$
contradicting  ($\star$).

To prove 4) consider the  plane $\pi$ osculating  $\ga_0$
at a point $p$. Assume that $\ga_0$ intersects $\pi$ at some point $q\neq p$.
Then $\pi$ is tangent to $\ga_0$ at $q$ (by item 1)),
$\ga_0(p)\neq \ga_0(q)$ (by item 2)) and the line $l$ passing through
the points $\ga_0(p),\ga_0(q)$ is not tangent to $\ga_0$ at these points
(by item 3)). There are two possible cases to consider:

i) the  plane osculating $\ga_0$ at the point $q$ is transversal to
the line $l$. In this case for any sufficiently small $t$
the  plane osculating the curve $\ga_t$ at some point $\tilde p$ close
to $p$ is tangent to $\ga_t$ at some point $\tilde q$ close to $q$.
It is clear, that the line passing through the points
$\ga_t(\tilde p),\ga_t(\tilde q)$ is not tangent to
$\ga_t$ at these points and intersects $D(\ga_t)$ with the
multiplicity  $> 4$.

ii) the plane $\pi$  osculates $\ga_0$ at $q$. Hence,
for any sufficiently small $t$ the  plane osculating $\ga_t$
at any point $\tilde p$ close to $p$ intersects $\ga_t$ at
some point $\tilde q$ close to $q$. As before the line
passing through the points $\ga_t(\tilde p),\ga_t(\tilde q)$
is not tangent to $\ga_t$ at these points and intersects
$D(\ga_t)$ with  multiplicity  $> 4$.

Thus we get a contradiction with  ($\star$) in both
    cases. \qed
\enddemo

As above denote by $L$ a line in $\bP^3$ which is not tangent to $\ga_0$
and intersects $D_0$ with total multiplicity $> 4$. The next lemma
describes restrictions on $L$ in the case when $\ga_{0}$ is locally convex
and satisfies ($\star$) and ($\star\star$).

\medskip
{\smc 2.2. Lemma.} {\it Under the above assumptions the line $L$ must
be tangent to every smooth irreducible
component of the surface $D_0$ at each intersection point of $L\cap D_0$.}

\demo{Proof} We will consequently exclude the possibility of having a
given number of transversal intersections. Namely, the line $L$ 
cannot  intersect
$D_0$ transversally at more than 4 smooth points since otherwise $L$ will
intersect $D(\ga_t)$ at least 5 times for all sufficiently
small $t$.

Let us show that $L$ cannot  intersect $D_0$ transversally
at 3 or 4 smooth points.
Indeed, in this case $L$ intersects $D_0$ at some point $P$ with the
local multiplicity greater than 1. If $P$ belongs to $\ga_0$
then  for any small $t$ any line close to $L$ and intersecting
the curve $\ga_t$ at a point close to $P$ will intersect $D(\ga_t)$
with  multiplicity $> 4$. If $L$ is tangent to a
smooth local irreducible component of the surface $D_0$ at the point $P$ then
again  for any small $t$ any line  close to $L$ and  tangent
to a smooth local irreducible
component of the surface $D(\ga_t)$ at a point near $P$ will
intersect $D(\ga_t)$ with multiplicity $> 4$.

Finally, suppose that $L$ intersects $D_0$ transversally  at $1$ or $2$
smooth points. Then $L$ intersects two local irreducible components
$D_0^1,D_0^2$ of the surface $D_0$ at some points $P_1,P_2$ with the
local multiplicities at least 2. By Lemma 2.1 items 2) and 4), there are
two cases to consider:

a) $P_1,P_2$ belong to the curve $\ga_0$ and $P_1\neq P_2$. In this
case  for any small $t$ the line passing through any points
$\tilde P_1,\tilde P_2$ of the
curve $\ga_t$ close to $P_1,P_2$ resp. intersects $D(\ga_t)$
with  multiplicity $> 4.$

b) $P_1,P_2$ do not lie on the curve $\ga_0$.
In this case the line $L$ is the intersection  of the tangent
planes to smooth components $D_0^1,D_0^2$ at the points
$P_1,P_2$ (these planes are transversal). Hence, for any small $t$
there exists a line in $\bP^3$  close to $L$ and tangent to
two smooth irreducible components of the surface $D(\ga_t)$
at the points near $P_1,P_2$. This line intersects $D(\ga_t)$
with the multiplicity $> 4$.

Thus, $L$ cannot  intersect $D_0$ transversally
at its smooth points. \qed

\enddemo

\medskip
Lemmas 2.1 and 2.2 leave us with only two possible (and dual to each
other) positions of the line $L$ w.r.t $\ga_{0}$:

\medskip
(i) the line $L$ intersects  $\ga_0$ at least at 3
different points and does not lie in osculating planes to $\ga_0$;

(ii) the line $L$ lies in at least  3 different osculating planes to
the curve $\ga_0$ and does not intersect $\ga_0$.

\medskip
The next lemma rules out these remaining possibilities and accomplishes
the proof of Proposition 1.

\medskip
{\smc 2.3. Lemma.} {\it In both cases {\rm (i)} and {\rm (ii)} one can
find a  line close to $L$ which  intersects $D_0$ transversally
  at more than  $4$ points.}

\demo{Proof} Indeed, consider the central projection
$\varrho:\bP^3\to \bP^2$ from some point $O\in L$. (In each of the cases we
will choose the point $O$ differently.)

(i) Assume that the point $O$ does not belong to $\ga_0$. Then
  $\varrho (\ga_0)$ has at least 3 irreducible locally
convex components $c_1,c_2,c_3$ at the point $\varrho (L)$. Every
curve $c_i, i=1,2,3$ separates a small neighborhood of this
point into two open domains one of which (denoted by $U_i$) satisfies
the property that each point of $U_i$ belongs to two lines  tangent to the
curve $c_i$. The intersection $U=U_1\cap U_2\cap U_3$ is nonempty.
  (See Fig.1a below illustrating the most delicate case.)
Then  for any point $P\in U$ the line $\varrho^{-1}(P)$
intersects at least 6 local irreducible components of $D_0$ transversally.

(ii) Let $L$ belong to the  planes osculating $\ga_0$ at the
points $P_1,P_2,P_3$. Then the tangent lines to  $\ga_0$
at these points cannot intersect the line $L$ at the same point
since the family $\ga_t$ is taken generic. Take the point $O\in L$
which belongs only to one of the
tangent lines to $\ga_0$ at the points $P_1,P_2,P_3$. Assume, for
example, that $O$ belongs to the line tangent  to $\ga_0$ at $P_1$.
Then the projection $c_i$ of a germ of
$\ga_0$ at the point $P_i$ has a half-cubical cusp at the point
$\varrho (P_i)$ for $i=1$ and cubical inflections for $i=2,3$.

It is easy to see that the projections of the
planes osculating $\ga_0$ at the points $P_1,P_2,P_3$ separate a small
neighborhood of the point $\varrho (L)$ into six parts one of which
(denoted by $U$) satisfies the property that each point of $U$ belongs to
two  lines tangent to each curve $c_i,i=2,3$ as well as one  line tangent
to the curve $c_1$, see Fig.1b.
Then  for any point $P\in U$ the line $\varrho^{-1}(P)$  intersects
at least 6 irreducible components of the surface $D_0$ transversally.
  Thus both cases (i) and (ii) are impossible.
This contradicts our initial assumption about the local convexity of $\ga_{0}$
and accomplishes the proof of Proposition 1. \qed
\enddemo

Proposition 1 implies that if a connected component of the space of
locally  convex curves contains a curve such that the degree of its
developable hypersurface is 4,
then the degree of the developable hypersurface of every curve in this
component is 4. In the (unique up to the choice of orientation) component
containing convex curves the degree of the developable hypersurfaces equals
$4$ since this fact holds for the rational normal curve. Hence
Theorem A follows from the next statement.

\medskip
    \vskip 15pt
\centerline{\hbox{\epsfysize=3,5cm\epsfbox{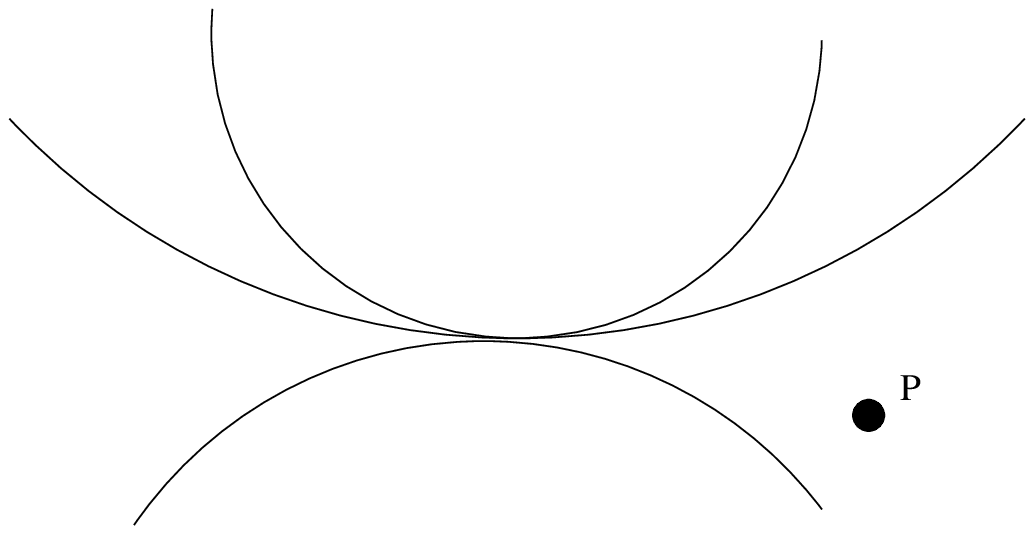}}
\hskip1cm\hbox{\epsfysize=4,5cm\epsfbox{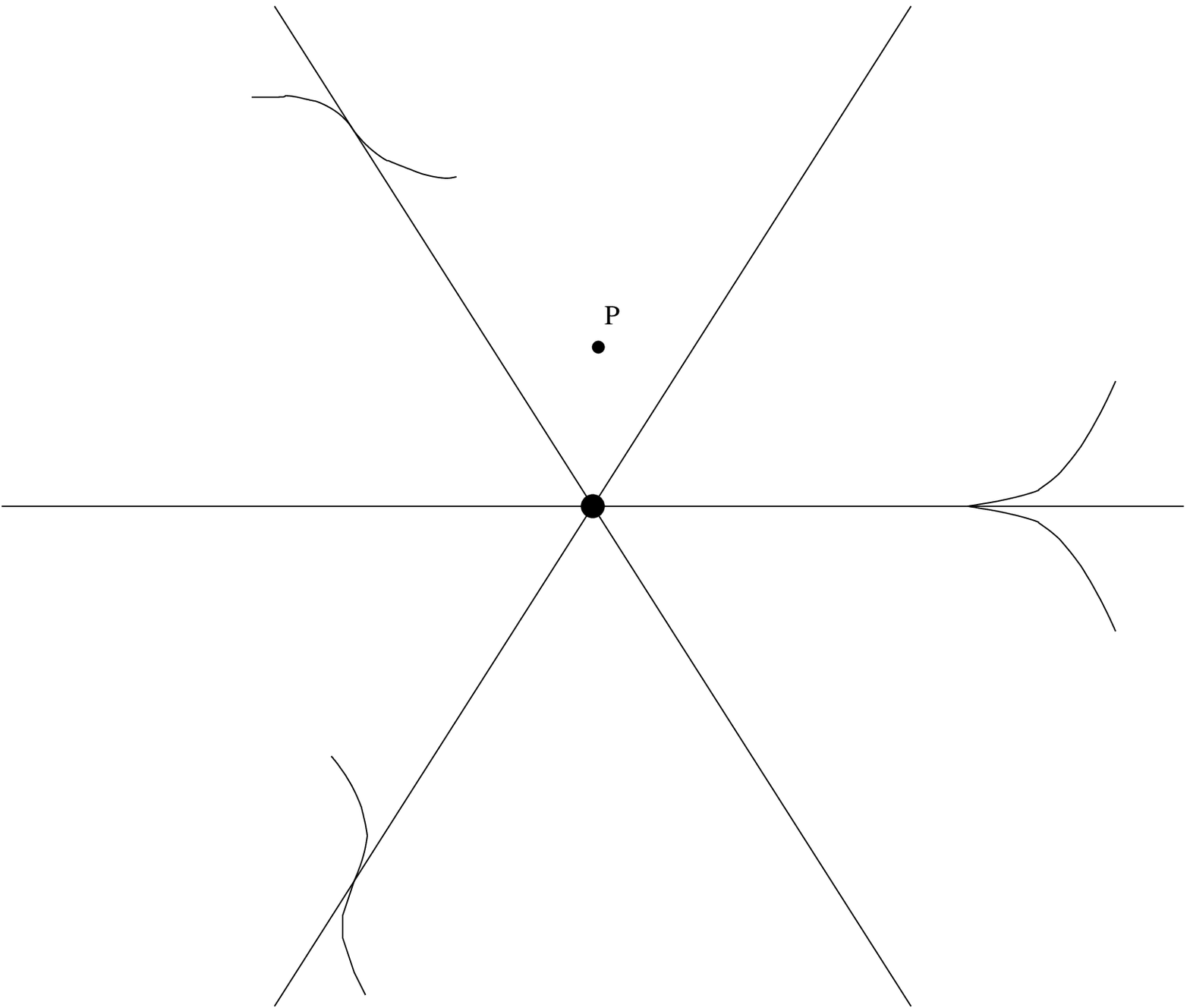}}
}
\midspace{0.1mm} \caption{Fig.~1. Illustrations to Lemma 2.3(i) and
2.3(ii).}

\medskip
{\smc Proposition 2.} {\it Every connected component of the space of
locally convex curves different from two components consisting of 
convex curves only
contains a curve with the degree of its developable hypersurface
  $> 4$.}

\demo {Proof}
     By \cite {MSh} there exist 5 connected components of the space of
     locally convex right projective curves in $\bP^3$. (The word
     `right' means that the Frenet frames of all these curves
     represent just one fixed `right' orientation of $\bP^3$.)
     Three of these components (including the component of convex curves)
     contain noncontractible curves and the other two contain 
contractible curves.
     The invariant distinguishing components is the mod 2 element of
     $\pi_{1}(SO_{4})$.

     Moreover, if $\xi$ is a contractible locally convex
     curve whose lift to $SO_{4}$ realizes a nontrivial element of
     $\pi_{1}(SO_{4})$ then we can choose the following  representatives of
     all five connected components. Namely, $\rho$, $\rho+\xi$ and $\rho+\xi^2$
     represent three components of noncontractible curves; $\xi$ and $\xi^2$
     represent two components of contractible curves. Here $\rho$ is the
     rational normal curve in $\bP^3$, the $'+'$-sign means that we
     append the second curve to the endpoint of the first curve with
     the same tangent direction.
     Finally, the expression $\xi^2$ denotes the curve $\xi$ traversed twice.
     (The lift of a
     locally convex right curve to
     $SO_{4}$ is obtained as follows, see also \S 3. A right  curve in
     $\bP^3$ first lifts  canonically to the sphere $S^3\subset \bR^4$ 
and then its
     extended Frenet frame lifts it to $SO_{4}$. One can lift all
     locally convex curves, both right and left to $O_{4}$.)

     Here and below by an {\it extended Frenet frame} of a projective curve
     $\ga: S^{1}\to \bP^3$ we mean the frame $(\tilde \ga, \tilde
     \ga', \tilde \ga^{\prime\prime}, \tilde 
\gamma^{\prime\prime\prime})$ of the
     canonical lift $\tilde \ga: S^1\to S^3\subset \bR^4$ of $\ga$ 
described above.

     Therefore if we can find an
     example of $\xi$ such that $D(\xi)$ has degree  $> 4$ Proposition
     2 will follow. An example of $\xi$ given in Fig. 2 below
     was obtained by a suitable modification of fig 2, p. 111  in \cite {To}.
     (The shown modification was required since the original picture 
in  \cite {To}
     which was supposed to present the
     velocity vector of an apropriate $\xi$ must necessarily
     contain spherical inflection points by the theorem of Fenchel
     on the existence of inflections on a noncontractible embedded
     curve on $\bP^2$, \cite{Fe}.)

     Let us now
     show that the developable hypersurface $D(\xi)$  has
     degree at least $6$. In fact, we show that already the piece
     $\tilde \xi$ of $\xi$ corresponding to the fragment of 
$\xi^\prime$ going around
     the north pole three times gives a contribution of at least $6$ to
     the degree. Indeed, consider the  $(x,y)$-plane $\Pi$. The
     velocity vector of $\tilde \xi$ is never vertical. Therefore the
     projection $\pi\xi$ of $\tilde \xi$ onto $\Pi$ along the $z$-axis   is
     smooth, convex and its velocity vector makes three complete turns.
     Let us choose some convex disc $\frak D$ containing $\pi\xi$. We
     check now that $\Pi\setminus D$ contains a point having at least
     $6$ different tangents to $\pi\xi$. Indeed, the union of all
     tangents to any  segment
     of a convex curve of $\bR^2$ whose tangent vector makes a halfturn covers
     the whole $\bR^2$ except for the strip between the
     parallel tangents at the endpoints. Therefore, the segment
     whose velocity vector makes three complete turns covers $\Pi$
     minus some strip at least $6$ times. Finally, choose the line $ll$
     parallel to the $z$-axis and passing through the chosen point
     on $\Pi$. We have just shown that the developable hypersurface
     $D(\xi)$
     of $\xi$ crosses $ll$ at least $6$ times.
     \qed
\enddemo

  \medskip
    \vskip 15pt
\centerline{\hbox{\epsfysize=4,5cm\epsfbox{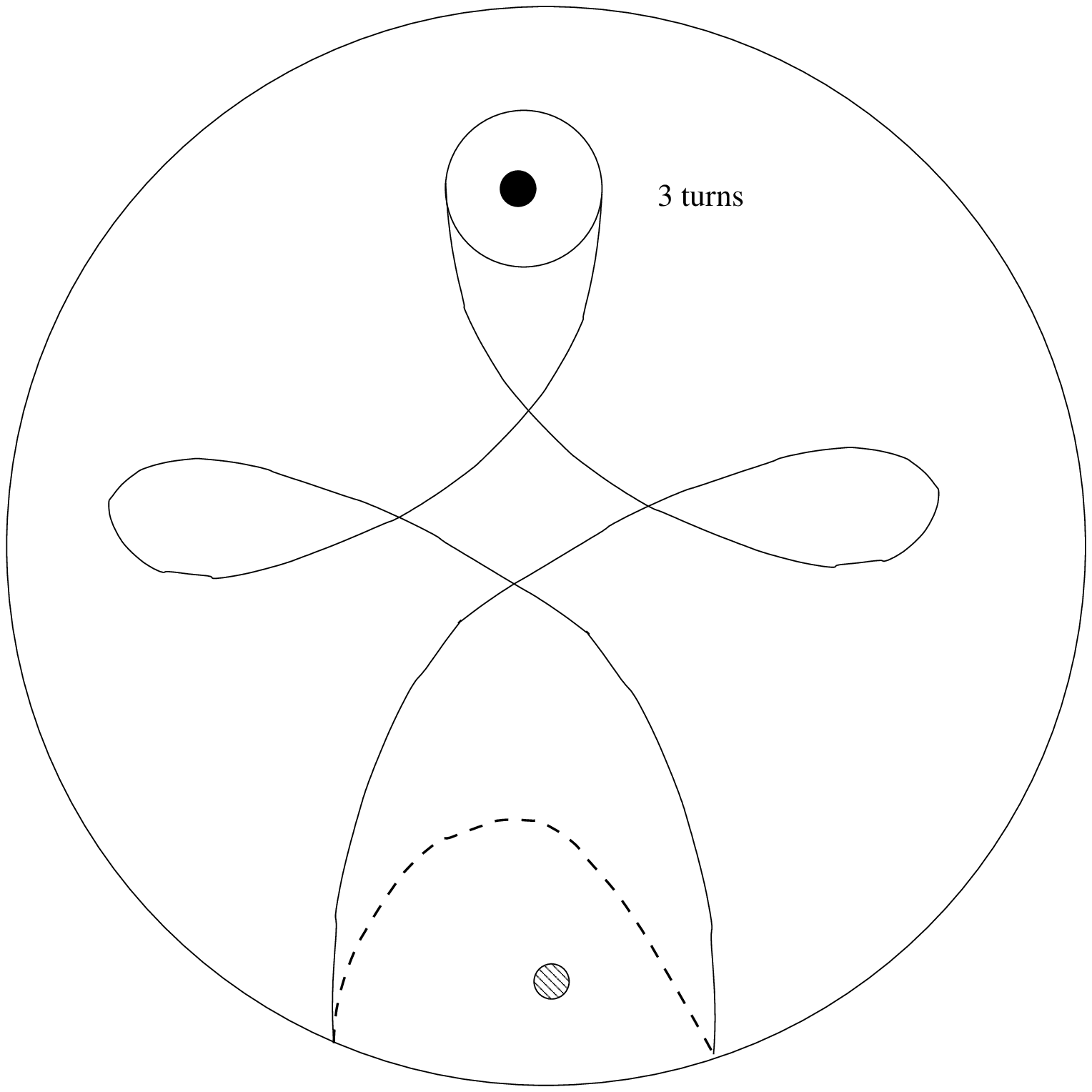}}}
\midspace{0.1mm} \caption{Fig.~2. Velocity vector of an appropriate
representative $\xi$. }

     \noindent
     {\bf Explanations to Fig.~2.}  This picture represents the
     graph of the velocity vector ${\xi^\prime}$ of $\xi$ parameterized by the
     arclength. Therefore $\vert {\xi^\prime} \vert =1$ implying that
     $\xi^\prime$ lies on $S^2\subset \bR^3$. (The black dots on the
     picture stand for
     the north resp. south poles.)
      The existence of a closed $\xi$ with a given $\xi^\prime$
     follows from a result of Fenchel's claiming that a closed
     nonparameterized curve in
     $\bR^n$ is realizable (with a suitable parametrization)
      as the velocity vector of another closed
     curve if and only if its convex hull contains the origin. The
     curve $\xi$ is contractible in $\bP^3$ (since it lies in $\bR^3$)
     and its lift to $SO_{3}$ by the usual Frenet frame realizes the
     nontrivial element in $\pi_{1}(SO_{3})$. This last fact follows
     from the following observation. We can deform $\xi^\prime$ in one
     hemisphere and the resulting curve makes an odd number of turns
     in this hemisphere. This implies using the
     triviality of the fibration $SO_{4} \to S^3$ restricted to
     a hemisphere that the lift of $\xi$ to $SO_{4}$ by the extended 
Frenet frame
     realizes a nontrivial element in $\pi_{1}(SO_{4})$ as well.

\heading \S 3. Proof of Theorem B. \endheading

In order to prove theorem B we first recall the relation between the
osculating flags to a convex curve and totally positive
uppertriangular matrices, see \cite {BSh}.

{\smc Some Notation.}
Let $F_{n+1}$ denote the space of complete (projective) flags in
$\bP^n$ or, equivalently, linear flags in $\bR^{n+1}$.
For a given locally convex curve $\ga: S^1\to \bP^n$ we define
  its {\it flag lift} $\ga_{F}: S^1\to F_{n+1}$ to be the curve of its
complete osculating flags. Two flags $f_{1}$ and $f_{2}$ are called
{\it transversal} if each pairs of subspaces (one from each flag) is
transversal. The set $Tn_{f}$ of flags nontransversal to a given
flag $f$ is called the {\it train } or the {\it flag hyperplane}
associated to $f$. Note that $Tn_{f}$ is the
union of all Schubert cells of positive codimension of the standard
Schubert cell decomposition of $F_{n}$ w.r.t $f$.

One has the following fundamental criterion of convexity going back to
G.~P\'olya, see \cite {Po}  and \cite {BSh}.
\medskip
{\smc 3.1. Theorem.} A curve $\ga: S^1\to \bP^n$ is convex if and
only if any  two distinct osculating flags are transversal.
\medskip

  A slight refinement of the above statement  allows us to characterize
  all configurations of complete flags which can be realized as
  osculating to some convex curve  in appropriate semigroup terms.
   We will use homogeneous coordinates and oriented complete flags.
   A complete linear flag in $\bR^{n+1}$ is called {\it oriented}
   if any of its subspaces is equipped with an orientation. The set of
   all complete oriented flags in $\bR^{n+1}$ is isomorphic to the
   group $O_{n+1}$ which covers $F_{n+1}$ with the discrete fiber of
   cardinality $2^{n+1}$ consisting of all choices of orientations for
   the subspaces. Notice that any projective curve $\ga: S^1\to \bP^n$ can be
canonically lifted to the curve $\tilde \ga$ on the sphere
$S^n\subset \bR^{n+1}$. (If $\ga$ is contractible then $\tilde \ga$ is
closed and if $\ga$ is noncontractible then $\tilde \ga$ ends at the
antipodal points. Note that convex curves in $\bP^{2n+1}$ are
noncontractible while convex curves in $\bP^{2n}$ are contractible.)
  Any locally convex $\ga: S^1\to \bP^n$ can therefore be lifted by
  using its extended Frenet frame in $\bR^{n+1}$
to the space $O_{n+1}$ of complete oriented flags. (The lift of a 
contractible curve
is a closed curve on $O_{n+1}$, the lift of a noncontractible is `antipodal'.)
As above we say
that two oriented flags $\tilde f_{1}$ and $\tilde f_{2}$ are {\it
transversal} if their underlying usual flags $f_{1}$ and $f_{2}$
(obtained by forgetting all orientations) are transversal. The
Schubert cell decomposition of $F_{n+1}$ w.r.t. some complete flag
lifts to $O_{n+1}$ and the lifted decomposition contains, for example,
$2^{n+1}$ open Schubert cells all projecting to
the only open Schubert cell in $F_{n+1}$.

Let us fix a basis $(e_{1},\ldots,e_{n+1})$ in $\bR^{n+1}$.
Let $f_{+}$ be the complete linear flag whose $i$-dimensional
subspaces are spanned by $e_{1},\ldots,e_{i}$ and $f_{-}$ be its
opposite flag whose $i$-dimensional subspaces are spanned by
$e_{n+1},\ldots,e_{n-i+2}$. Analogously, let $\tilde f_{+}$ be the complete
oriented linear flag whose $i$-dimensional
subspaces are spanned by $e_{1},\ldots,e_{i}$ with the orientation
given by $e_{1}\wedge \ldots \wedge e_{i}$ and $\tilde f_{-}$ be its
opposite flag, i.e. the flag whose $i$-dimensional subspaces are spanned by
$e_{n+1},\ldots,e_{n-i+2}$ with the orientation given by
$e_{n+1}\wedge \ldots \wedge e_{n-i+2}$.
  The open Schubert cell (affine chart) on
$F_{n+1}$ consisting
of all flags transversal to $f_{-}$ can be identified with the group $T$
of all uppertriangular $(n+1)\times(n+1)$-matrices (taken in the
chosen basis $(e_{1},\ldots,e_{n+1})$) whose diagonal entries equal 
$1$.  Analogously, the
open Schubert cell on $O_{n+1}$ consisting
of all oriented flags transversal to $f_{-}$ and containing $\tilde
f_{+}$ can be identified with the group $T$ as well. (As usual, the 
$i$-dimensional
subspace of the flag corresponding to some matrix is spanned by its $i$
first rows and its orientation is induced by the wedge of these rows.)

{\smc Definition.} A matrix $M\in T$ is called uppertriangular {\it
totally positive} (resp. {\it totally nonnegative} ) if any of its
minors  which does not vanish due to uppertriangularity is positive
(resp.  any of its minors is nonnegative), see e.g. \cite {Ga}.
\medskip
{\smc Remark.} The set $T^+\subset T$ of all totally
positive matrices is a semigroup and is contractible as a topological
space.
\medskip
  {\smc Definition.} A sequence $\{\tilde f_{1}, \tilde
  f_{2},\ldots , \tilde f_{r}\}$ of matrices in $T^+\subset T$ is
  called {\it totally positive} if each $\tilde f_{i+1}$ lies in the
  cone $\tilde f_{i} \circ T^+$, i.e. there exists a totally positive matrix
  $M_{i,i+1}$ such that $\tilde f_{i+1}=\tilde f_{i} M_{i,i+1}$. If
  $T^+$ is identified  with is the open cell in $O_{n+1}$ w.r.t.
   some basis $(e_{1},\ldots,e_{n+1})$ in $\bR^{n+1}$ then
  we say that the sequence  of oriented flags $\{\tilde f_{1}, \tilde
  f_{2},\ldots , \tilde f_{r}\}$ is {\it totally positive w.r.t. this basis}.

Let us choose an orientation of $S^1$ and a reference point $O\in S^1$.
This choice identifies $S^1\setminus O$  with the interval $[0,2\pi)$.

\medskip
{\smc 3.2. Theorem, see \cite {BSh}.}
For a given convex curve $\ga: S^1\to \bP^n$ and a given sequence of
points $\{t_{1}<t_{2}<\ldots < t_{r}\} \in [0,2\pi)=S^1\setminus O$ 
there exists a
basis in $\bR^{n+1}$ such that the sequence $\{\tilde f_{1},\ldots,
\tilde f_{r}\}$ of extended osculating
flags to $\tilde \ga$ at the points  $\{t_{1}<t_{2}<\ldots < t_{r}\}$
  is totally positive  w.r.t. the above mentioned basis. Conversely,
for any given sequence of totally matrices
there exists a convex curve $\ga: S^1\to \bP^n$ and a
set of points $\{t_{1}<t_{2}<\ldots < t_{r}\}\in [0,2\pi)$ such that
  this sequence coincides with the sequence of extended osculating
  oriented flags to this curve at these points w.r.t. an appropriately
  chosen basis.

\demo{Proof} Let us sketch the proof, see details in \cite {BSh}. The
  lifts of locally convex curves to $O_{n+1}$ obey the
following remarkable {\it Cartan} distribution of cones on $O_{n+1}$.
Let $\tilde f=(\bold l_1 \subset \dots \subset \bold l_{n+1})$ be a complete
oriented flag on $\bR^{n+1}$, where $\bold l_{i}$ denotes its
$i$-dimensional oriented subspace. We define in $O_{n+1}$
the set of
$n$ circles $\{c_1(\tilde f) , \dots , c_{n}(\tilde f)\}$ passing through
$\tilde f$ and given
by the
relation
$$ c_i(\tilde f) =\{ \bold l_1 \subset \dots \subset \bold l_{i-1} \subset L_i
\subset \bold l_{i+1} \subset \dots \subset
\bold l_{n+1} \}\qquad i= 1,\ldots,n     $$
where $L_i$ runs over the set of all $i$-dimensional subspaces 
satisfying the above
inclusions and having the appropriate orientation, i.e. such that
$\tilde f$ itself belongs to $c_{i}(\tilde f)$.

The tangent lines to $c_1(\tilde  f) , \dots , c_{n}(\tilde f)$ at
$\tilde f$ are linearly  independent. They also have prescribed
orientations meaning that the orientation of $\bold l_{i}$ appended by
the velocity vector of the rotation of $L_{i}$ around  $\bold
l_{i-1}$ gives the orientation of $\bold l_{i+1}$. Let $\tenit 
c_{i}^+(\tilde f)$
denote the chosen  tangent halfline to $c_{i}(\tilde f)$.
Take the open $n$-dimensional orthant $ C_{\tilde f}
\in {TO}_{n+1}$ spanned by all $\tenit c_{i}^+(\tilde f),\; i=1,\ldots ,n$.
  The distribution $\frak C=\bigcup_{\tilde f \in O_{n+1}} C_{\tilde f} $
is called the {\it Cartan distribution} on the space $O_{n+1}$.
  In Lie-theoretic terms this distribution can be described  as
  follows. $\frak C$ coincides with
  the left-invariant distribution on the group $O_{n+1}$ which is
  generated by the cone of skew-symmetric matrices having all
  vanishing entries except those on the 1st upperdiagonal which
  are positive and, respectively, those on the 1st lowerdiagonal which
  are negative.
  \medskip
  {\smc 3.2.1. Lemma, see \cite {BSh}.} The flag lift of any locally
  convex curve to $O_{n+1}$ is tangent to $\frak C$.
  \medskip

  It is fairly easy to see that the
image of this Cartan distribution restricted to the group $T$ 
(considered as the open
cell in $O_{n+1}$) coincides with the left-invariant distribution generated
by the cone of all positive linear combinations of the entries on the
1st upperdiagonal in the  Lie algebra $\frak T$ of $T$ consisting of 
all nilpotent uppertriangular
matrices. We call this distribution the {\it Cartan distribution} on
$T$ and by abuse of notation we denote it by $\frak C$ as well.
  The following result is classical and crucial in understanding of Theorem 3.2.
\medskip
{\smc 3.2.2. Loewner-Whitney theorem}, see e.g. Th. 1.1. in \cite 
{FZ}. The Cartan distribution
$\frak C$ on $T$ has
the property that the accessibility domain of any matrix $m\in T$ (i.e. the
set of all end points of the curves starting at $m$ and tangent to $\frak
C$) coincides with $m\circ T^{+}$.
\medskip The idea of the proof is as follows. Let us fix any
irreducible decomposition $\bar i=(i_{1},\ldots, i_{\binom {n}{2}})$ 
of the longest permutation
$(n,n-1,n-2,\ldots,1)$. Define the map $\Psi_{\bar i}: 
(\bR^+)^{\binom {n}{2}}\to T$ which
sends a $\binom {n}{2}$-tuple of positive numbers
$(\tau_{1},\ldots,\tau_{\binom {n}{2}})$ to the product of
$T_{\tau_{1},i_{1}} \circ T_{\tau_{1},i_{1}} \circ \ldots \circ T_{\tau_{\binom
{n}{2}},i_{\binom{n}{2}}}$ where $T_{\tau_{k},i_{k}}$ is the
uppertriangular unipotent matrix whose only nonzero element except the
main diagonal equals ${\tau_{k}}$ and stands at the $i_{k}$-th row
and $i_{k+1}$-st column. One can show that the image of $\Psi_{\bar
i}$ coincides with $T^{+}$ with some subset of  real codimension at
  least $2$ removed. Thus the closure of the image equals the set of all
totally nonnegative matrices. This guarantees that the accessibility domain
of the identity matrix is the whole $T^{+}$ which together with the 
left-invariance of our
distribution implies the required result.
\medskip
This statement together with Theorem 3.1 and the interpretation of
$T$ as an open Schubert cell in $O_{n+1}$ prove Theorem 3.2. \qed

\enddemo

\medskip

The next two simple lemmas can be found in e.g. \cite {VD}.
\medskip
{\smc 3.3. Lemma.} For any triple of pairwise nonintersecting lines
$l_{1},l_{2},l_{3}$ in $\bP^3$ there exists a unique hyperboloid
$H$ of one sheet containing $l_{1},l_{2},l_{3}$.
\medskip

Take any  $4$-tuple $L$ of pairwise nonintersecting lines
$L=(l_{1},l_{2},l_{3},l_{4})$ in $\bP^3$.
\medskip
{\smc 3.4. Lemma.} There exist two real lines intersecting each line in $L$
if and only if $l_{4}$ intersects  in two distinct points the
hyperboloid $H$ containing $l_{1},l_{2},l_{3}$.
(This property is independent on the particular choice of a triple in $L$. )
\medskip

{\smc 3.5. Main example.} The following $4$-tuple of lines in
$\bP^3$ (given in homogeneous coordinates in $\bR^4$) satisfies
the conditions of Theorem 3.2 and settles Theorem B.

$$l_{1}=\pmatrix 1&0&0&0\\
                   0&1&0&0 \endpmatrix, \quad
    l_{2}=\pmatrix 1&1&1/2&1/6\\
                   0&1&1&1/2 \endpmatrix, $$

$$  l_{3}=\pmatrix 1&101/100& 401/200& 2503/1500\\
                   0&1&5/2&5/2 \endpmatrix, \quad
   l_{4}=\pmatrix 0&0& 1&0\\
                   0&0&0&1 \endpmatrix. $$

   The totally positive uppertriangular matrices $M_{i,i+1}$ sending $l_{i}$ to
   $l_{i+1}$ are as follows

   $$M_{1,2}=\pmatrix 1&1&1/2&1/6\\
                      0&1&1&1/2\\
		    0&0&1&1\\
		    0&0&0&1\endpmatrix,\quad
    M_{2,3}=\pmatrix 1&1/100&1/200&1/500\\
                      0&1&3/2&1\\
		    0&0&1&1\\
		    0&0&0&1\endpmatrix. $$
  Finally, $l_{4}=\lim_{t\to +\infty}l_{3} M_{3,4}(t), $
   where
  $$ M_{3,4}(t)= \pmatrix 1&t&t^2/2&t^3/6\\
                      0&1&t&t^2/2\\
		    0&0&1&t\\
		    0&0&0&1\endpmatrix . $$
\medskip

\demo{Proof} Verified by direct calculation. \qed
\enddemo
\medskip

{\smc 3.6. Lemma.} The line $l_{3}$ does not intersect the hyperboloid
$H\subset \bP^3$ containing $l_{1},l_{2},l_{4}$.

\demo{Proof} Indeed, the hyperboloid $H$ is given by the equation
$\Phi=2 x_{1}x_{3}-3 x_{1} x_{4} - 3 x_{2} x_{3} +6 x_{2} x_{4}$.
  If we apply to the line $l_{2}$ an uppertriangular matrix
$$M=\pmatrix 1&a&b&c\\
               0&1&d&e\\
	      0&0&1&f \\
	      0&0&0&1\endpmatrix$$
we get the line $M(l_{2})$ given by $\pmatrix 1&a+1&b+d+1/2&c+e+f/2+1/6\\
                                           0 & 1 & d+1 &
                                           e+f+1/2\endpmatrix. $

Take the affine line in $\bR^4$ $$ll=\cases x_{1}=1\\
                                  x_{2}=1+a+u\\
				x_{3}=1/2+b+d+u(d+1)\\
				x_{4}=1/6+c+e+f/2+u(1/2+f+e)\endcases$$
belonging to $l_{3}$ (considered as a $2$-dimensional subspace in
$\bR^4$). The restriction of the defining equation of $H$ to $ll$
(calculated using Mathematica)  is given by
$$\align \Phi(ll)=& u^2(-3 d + 6 e + 6 f)+u(-3 b + 6 c - 4 d - 3 a d 
+ 9 e + 6 a e + 6 f +
6 a f)  - \\ & a/2 - b - 3 a b + 3 c + 6 a c  - d - 3 a d + 3 e +
     6 a e + 3 f/2 + 3 a f .\endalign$$

The discriminant of this quadratic in $u$ equation  equals
$$\align Dsc(ll)=& 9 b^2 - 36 b c + 36 c^2 - 6 ad + 12 b d - 18 a bd -
     12 c d + 36 a cd +  4 d^2 - 12ad^2 + 9 a^2 d^2 +
     12 a e -\\ & 30 b e + 36 a b e + 36 c e - 72 a c e - 12 d e +
     42 a d e - 36 a^2 de + 9 e^2 - 36 a e^2 + 36 a^2 e^2 +
     12 a f - \\ & 12 b f + 36a b f - 72 a c f - 6 d f + 24 ad f -
     36 a^2 d f - 36 ae f + 72 a^2 e f + 36 a^2 f^2. \endalign$$

Observe that if $a=b=c=0$ then $Dsc_{ll}=4 d^2 - 12 d e + 9 e^2 -
6 d f=(2d-3e)^2-6df.$ If one additionally assumes that $d>0, e>0,
f>0$ and $df>e$ then the corresponding matrix $M$ is totally nonnegative. One
can easily adjust the values of $d,e,f$ to get the negative value of
$Dsc(ll)$. For example, $d=3/2, e=f=1$ gives $Dsc(ll)=-9.$
Since the space of all totally nonnegative matrices is the closure of
the space of all totally positive matrices (see \cite {Wh}) one can 
find a small
deformation of the later matrix which is totally positive but
$Dsc(ll)$ is still negative. The matrix $M_{2,3}$ above is totally
positive with $l_{3}=l_{2}M_{2,3}$ such that
$Dsc(ll)=-2231979/250000\approx -8.92792.$ \qed

\enddemo

To complete the proof of Theorem B notice that by Theorem 3.2
there exists a convex curve whose tangent lines are
$l_{1},l_{2},l_{3},l_{4}$ and thus by lemma 3.4 there are no real
projective lines intersecting all four of them. In fact, in order to
get the exact set-up of Theorem 3.2 we need to take instead of the
line $l_{4}$ the line $l_{3}M_{3,4}(t)$ for some sufficiently large
$t$, see the above notation.
\qed

\heading \S 4. Some Grassmann geometry, etc. \endheading

\medskip
{\smc $1$st conjecture.} Let as above $G_{{k+1,n+1}}$ denote the 
usual Grassmannian of $(k+1)$-dimensional
real linear subspaces in $\bR^{n+1}$ (or equivalently, 
$k$-dimensional projective
subspaces in $\bP^{n}$).

{\smc Definition.} Given an $(n-k)$-dimensional subspace $L\subset \bR^{n+1}$
we define {\it the Grassmann hyperplane $H_{L}\subset G_{k+1,n+1}$ 
associated to
$L$} to be the set of all
$(k+1)$-dimensional subspaces in $\bR^{n+1}$ nontransversal to $L$.
\medskip
{\smc Definition.} A smooth closed curve $\ga: S^1\to G_{k+1,n+1}$ is called
{\it locally
Grassmann convex} if the local multiplicity of intersection of $\Ga$
with any Grassmann
hyperplane $H_L\subset G_{k+1,n+1}$ does not exceed $(k+1)(n-k)=\dim
G_{k+1,n+1}$ and
{\it globally Grassmann convex} or simply {\it Grassmann convex} if
the above condition holds for the global multiplicity.

Given a locally convex curve $\ga:S^1\to \bP^n$ and $1\le k\le n-1$
we can
define {\it its $k$-th Grassmann lift }
       $\ga^k_{G}: S^1\to G_{k+1,n+1}$ formed by the osculating
$k$-dimensional projective
subspaces to the initial $\ga$. (The curves $\ga^k_{G}$ are well-defined
for any $k=1,...,n-1$
due to the local convexity of $\ga$.)

The following is the reformulation of the first basic conjecture from
the introduction.
\medskip
{\bf Conjecture on the $k$-th developable restated.} For any convex $\ga:S^1\to
\bP^n$ its $k$-th Grassmann lift $\ga^k_{G}$ is Grassmann convex.
\medskip

   The scheme of the proof of Theorem A in \S 2 can be generalized to
   higher dimensions but the number of cases increases drastically.
   One apparently needs a different inductive argument.

   \medskip
{\smc $2$nd conjecture.}  Recall that an enumerative
problem in complex algebraic geometry is called {\it totally
real\, } if all objects involved can be realized in the real subspace of the
complex space. There is a large number of enumerative problems known to be
      totally real. On the other
hand, there exists an impressive list of those which fail to be
totally real, see  \cite {So3} and references therein.

The following question was raised by W.~Fulton in \cite {F}.

{\smc Problem.} Do there exist configurations of $(k+1)(n-k)$-tuples of
$k$-dimensional subspaces in $\bP^n\subset \bcP^n$ such that all $\sharp_{k,n}$
possible complex $(n-k-1)$-dimensional {\it complex} subspaces in
$\bcP^n$ intersecting each subspace in the $(k+1)(n-k)$-tuple are  real?
(Such configurations of subspaces are called {\it totally real}.)
\medskip
{\bf Conjecture on total reality restated.}  For the rational normal curve
$\rho_{n}: S^1\to \bP^{n}$ the configuration of
any   its $(k+1)(n-k)$ pairwise distinct projective $k$-dimensional
osculating subspaces is totally real.

Identifying the space $\bP^n$ with the space of all polynomials in
one variable of degree at most $n$ (considered up to a constant 
factor) and the rational normal curve with
the family $(x-\alpha)^n$ one gets the reformulation of the conjecture
on total reality in terms of rational curves given in \S 1. The most
  promising approach to this conjecture seems to be the study of the
geometry of the Wronski map from Grassmannians to projective spaces,
see e.g. \cite {EG2}. Namely, consider the map $W_{k,n}:
G_{k+1,n+1}\to \bP^{(k+1)(n-k)}$ associating to a
$(k+1)$-dimensional linear subspace in the space of polynomials
its Wronskian, i.e. the determinant of the $(k+1)\times (k+1)$ matrix
$$\pmatrix p_{1}(x) & p_{2}(x) &\ldots & p_{k+1}(x)\\
            p'_{1}(x) & p'_{2}(x) &\ldots & p'_{k+1}(x) \\
	    p''_{1}(x) & p''_{2}(x) &\ldots & p''_{k+1}(x)\\
	    \vdots & \vdots & \ddots &\vdots \\
	     p_{1}^{(k)}(x) & p_{2}^{(k)}(x) &\ldots & p_{k+1}^{(k)}(x)
	     \endpmatrix $$
	     where $p_{1}(x),\ldots, p_{k+1}(x)$ is some basis of the chosen
	     subspace.

   Note that any change of a basis of a given subspace
   results in multiplying the Wronskian by a constant.  Wronskians are
   polynomials of degree at most $(k+1)(n-k)$.
   (The Wronski map is much more natural over $\bC$ but
	     we are doing real algebraic geometry at the moment.)
	     One can now define two
   hypersurfaces in the image space $\bP^{(k+1)(n-k)}$ (identified
   with the space of polynomials of degree at most $(k+1)(n-k)$.)
    The first hypersurface $\Cal D\subset \bP^{(k+1)(n-k)}$ is the
    usual discriminant consisting of all polynomials with multiple
    zeros. The second hypersurface $\Cal {J} \subset \Bbb P^{(k+1)(n-k)}$
    is the image of the hypersurface in $G_{k+1,n+1}$ where the
    Jacobian of $W_{k,n}$ vanishes, i.e. where $W_{k,n}$ is not a local
    diffeomorphism.  If one could settle the following conjecture then 
together with the main
    result of \cite {So2} this would prove the conjecture on total reality.
    \medskip

    {\bf $\Cal D -\Cal J$-conjecture.} The hypersurface $\Cal J$ does
    not intersect the domain of all strictly hyperbolic polynomials,
    i.e.  polynomials whose zeros are all real and
    distinct.

    {\smc Problem.} Enumerate irreducible components of $\Cal D \cap
    \Cal J$.

    \medskip
{\smc On projections of convex curves.}
The following question seems to be of importance in computer graphics 
applications.

{\smc Problem.} Given a real plane closed smooth curve $\ga$ such that the
supremum of the total number of its intersection points with affine
lines is finite represent $\ga$ as a projection of a convex curve lying
in the projective space of some dimension.

At the moment it is not even proved (although very plausible)
that such a representation always exists. Assuming that it exists one
wants to determine  the minimal dimension  such that a given plane curve can
be obtained
projecting a convex curve lying in  the projective space of this dimension.
The obvious guess
(which the authors encounted in the applied literature) that this dimension
coincides with the above mentioned supremum of total number of
intersections is definitely wrong. We
finish the article with the following statement.

{\bf Conjecture.} A closed real smooth curve $\ga$ in $\Bbb P^2$ can be
obtained as a projection of a convex curve in $\Bbb P^3$ if and only
if the following holds:

\noindent
i) the maximal total intersection of $\ga$ with any line in $\Bbb
P^2$ equals 3;

\noindent
ii) the maximal number of inflection points of $\ga$ equals $3$;

\noindent
iii) the maximal number of tangent lines to $\ga$ through any point on
$\Bbb P^2$ equals $4$.

\enddocument